\documentclass[preprint,10pt]{elsarticle}
\usepackage{amsfonts,amsmath,amssymb,epsfig,graphicx,algorithm,subfig,cancel}
\usepackage{booktabs, wrapfig}
\usepackage{color,soul}

\usepackage[export]{adjustbox}

    \newtheorem{thm}{Theorem}
    \newtheorem{lem}[thm]{Lemma}
    \newdefinition{rmk}{Remark}
    \newdefinition{defi}{Definition}
    \newproof{pf}{Proof}

\makeatletter
    \newcommand{\Rmn}[1]{\it{\expandafter\@slowromancap\romannumeral #1@}}
\makeatother

\journal{Journal of Computational Physics}

\newcommand{\J}{\langle J \rangle}
\newcommand{\w}{{\bf w}}
\newcommand{\vv}{{\check v}}
\newcommand{\ww}{{\check w}}

\begin{document}
\begin{frontmatter}

\title{Simplified Least Squares Shadowing sensitivity analysis for
chaotic ODEs and PDEs}

\author[mit]{Mario Chater\corref{cor}}
\ead{chaterm@mit.edu}
\author[mit]{Angxiu Ni}
\ead{niangxiu@mit.edu}
\author[mit]{Qiqi Wang}
\ead{qiqi@mit.edu}

\cortext[cor]{Corresponding author.}
\address[mit]{Aeronautics and Astronautics,
MIT, 77 Massachusetts Ave, Cambridge, MA 02139, USA}

\begin{abstract}
This paper develops a variant of the Least Squares Shadowing (LSS)
method, which has successfully computed the derivative for several
chaotic ODEs and PDEs.  The development in this paper
aims to simplify Least Squares Shadowing method
by improving how time dilation is treated.
Instead of adding an explicit time dilation term as in the original
method, the new variant uses windowing, which can be more efficient and simpler to
implement, especially for PDEs.
\end{abstract}


\end{frontmatter}

\renewcommand\u{\mathbf{u}}
\newcommand\R{\mathbf{R}}

\section{Introduction}
\label{s:intro}

In many engineering and non-engineering applications, quantities of interest are time averages of some specific instantaneous output and this output may be affected by parameters of the physical system/model. For example, when evaluating the performance of a wing, the time averaged drag it produces is of great importance for aerodynamicists. Any shape parameter of this wing (width, length, curvature...) clearly has an impact on the produced drag and could be consequently considered a parameter of the system. Computing the derivative of the quantity of interest (time-averaged output) with respect to the system's parameters is crucial in:
\begin{itemize}
\item \textbf{Numerical optimization}: the derivative feeds a gradient-based algorithm which is then used to optimize the systems parameters \cite{nocedal},\cite{designopti},\cite{airfoil}
\item \textbf{Uncertainty quantification}: the derivative itself gives a useful assessment of the sensitivity and/or uncertainty of the system with respect to its parameters \cite{rcralph},\cite{uncertainty}
\end{itemize}

For chaotic systems such as those encountered in turbulent fluid flows, the so-called \textit{butterfly effect} makes conventional tangent/adjoint algorithms fail when computing the derivative of the quantity of interest. Indeed, the computed derivatives of infinite-time averaged quantities are, in general, orders of magnitude bigger than the real ones \cite{butterfly}. Some methods have been proposed to overcome this problem \cite{ensembleadjoint},\cite{abramov},\cite{cooper}. \textit{Least Squares Shadowing} (LSS) is a relatively simple and efficient algorithm that evaluates the desired derivative by solving a constrained least squares problem \cite{explosion}. In the original version, corrective time dilation factors have to be introduced to ensure its convergence but these same terms lead to a heavier and more expensive system of equations to solve. Indeed, the presence of time dilation terms prohibits the use of sparse solvers in a crucial step of the LSS algorithm which consists of inverting a matrix. Since this matrix scales up with the dimension of the problem, it is necessary to address this issue if we want to improve our algorithm for very large, industrial-scale problems \cite{fluidexample}.
In this paper, we introduce a modified version of LSS in which the time dilation factors are replaced by a simple \textit{windowing} procedure: the result of the non corrected system of equations is multiplied by a specific windowing function. This procedure has been used for systems that have a periodic behavior \cite{periodic}, but we will show that the fundamental idea is still valid in our framework and can be transposed to chaotic dynamical systems. Since this \textit{"windowing"} operation is almost costless, the new algorithm is much more efficient than standard LSS, specially for high-dimensional systems.\\

The paper is organized as follows: in the next section, we introduce the new algorithm in both its tangent and adjoint versions. Then, in the third section, the underlying theory of the algorithm is presented and we prove that the result given by simplified LSS converges to the one given by the original algorithm when integration time goes to infinity. Finally, the last section is dedicated to a numerical application of the new algorithm using the Lorenz $63$ test case dynamical system.

\section{A summary of the tangent and adjoint versions of simplified Least Squares Shadowing}
\label{s:algo}

Consider a family of
ergodic dynamical systems satisfying a differential
equation parameterized by $s$ :
\begin{equation} \label{ode}
\frac{du}{dt} = f(u,s)
\end{equation}
where $f$ is smooth with respect to both $u$ and $s$.  
Ergodicity means that an infinite-time averaged quantity
\[ \J(s) := \lim_{T\to\infty} \frac1T \int_0^T J(u,s) dt \]
where $J$ is a measurable function and $u=u(t;s)$ is a trajectory for parameter
$s$, does not depend on the initial condition of the trajectory $u(0;s)$.
The derivative of such infinite-time averaged quantity can be approximated
by the following three steps:
\begin{enumerate}
    \item Integrate equation (\ref{ode}) until the initial transient behavior has
        passed.  Then further integrate the equation over a period of
        $T$, which should be a multiple of the longest time scales in its
        solution.  Store the trajectory of the second time integration
        as $u(t),0\le t\le T$.
    \item Find the function $\vv(t),0\le t\le T$ which minimizes $\int_0^T |\vv(t)|^2dt$ and 
        satisfies the linearized equation
        \begin{equation}\frac{d\vv}{dt} = f_u \vv + f_s,
        \end{equation}
        where the time-dependent Jacobians $f_u:=\frac{\partial f}{\partial u}$ and
        $f_s:=\frac{\partial f}{\partial s}$ are based on the trajectory
        obtained in the previous step.  This $\vv(t)$ can be found by
        solving a system of linear equations derived from the KKT conditions of the
        constrained least squares problem :
        \begin{equation} \label{kkttan} \begin{cases}
            \frac{d\vv}{dt} = f_u \vv + f_s & \\
            \frac{d\ww}{dt} = -f_u^T \ww + \vv&  \\
            \ww(0) = \ww(T) = 0             & 
        \end{cases}\end{equation}
        Note that this is the Least Squares Shadowing problem in
        previous literature \cite{explosion}, but without a time-dilation term.
    \item Approximate the derivative of $\J(s)$ by computing a windowed
        time-average :
        \begin{equation} \label{objlsstan}
            \frac{d\J}{ds}\approx \frac1T\int_0^T
            \w(\tfrac{t}{T}) \left(J_u \vv + J_s \right) dt
        \end{equation}
        where $\w$, the window function, is a scalar function in $[0,1]$ satisfying 
        \begin{enumerate}
            \item $\w$ continuously differentiable,
            \item $\w(0)=\w(1)=0$,
            \item $\int_0^1 \w(r) dr = 1$.
        \end{enumerate}
        An example satisfying all three criteria is $\w(r) = 1 - \cos 2 \pi r$.
\end{enumerate}
This resulting derivative approximation converges to the true
derivative as $T\to\infty$, as mathematically derived in Section \ref{proof}
and under the same assumptions\footnote{The essential assumption is for the dynamical system to be uniformly hyperbolic.} as the original Least Squares Shadowing method \cite{proof}.\\

In the algorithm above, the cost of solving the constrained
least squares problem (\ref{kkttan}) scales with the dimension of $s$, and is
independent from the dimension of $J$.  Such an algorithm, which favors
a low-dimensional parameter $s$ and a high-dimensional quantity of
interest $J$, is called the tangent (or forward) version of simplified LSS.\\
A corresponding adjoint (or backward) version can be derived,
whose computation cost favors a high-dimensional $s$ and a
low-dimensional $J$. This adjoint algorithm, also consists of three steps :
\begin{enumerate}
    \item Obtain a trajectory $u(t)$, $0\le t\le T$ in the same way as Step
        1 of the previous algorithm.
    \item Solve the system of linear equations
        \begin{equation} \label{kktadj} \begin{cases}
            \frac{d\hat{\ww}}{dt} = f_u \hat{\ww} & \\
            \frac{d\hat{\vv}}{dt} = -f_u^T \hat{\vv} - \hat{\ww} -
            \w(\tfrac{t}{T}) J_u^T &  \\
            \hat{\vv}(0) = \hat{\vv}(T) = 0             & 
        \end{cases}\end{equation}
        where $\w(r)$, $0\le r\le 1$ is a scalar windowing function
        satisfying the criteria described in Step 3 of the previous
        algorithm.  This system of linear differential equations is
        the dual of the system in Step 2 of the previous
        algorithm, derived by combining it with Equation (\ref{objlsstan})
        and integrating by parts. A complete derivation of the adjoint version can be found in \ref{adjointderivation}. 
    \item Approximate the derivative of $\J(s)$ by the following
        equation, derived together with Equations (\ref{kktadj}) in Step 2 :
        \begin{equation}\label{obffunadj}
            \frac{d\J}{ds}\approx \frac1T\int_0^T \left(f_s^T \hat{\vv} + J_s \right) dt
        \end{equation}
\end{enumerate}
We can show that the adjoint version of LSS (equations \ref{kktadj},\ref{obffunadj}) produces the same estimation of $\frac{d\J}{ds}$ than the one given by the tangent version (equations \ref{kkttan}, \ref{objlsstan}) up to round-off errors \cite{explosion}. The approximated
$\frac{d\J}{ds}$ should be the same if the differential
equations are solved exactly in both algorithms, and the integrals are
evaluated exactly.  With such exact numerics, the error in the
approximation is solely due to the infeasibility of using an infinite $T$,
and should diminish as $T\to\infty$.
This error depends on the trajectory $u(t)$, but does not depend on
whether the tangent or adjoint algorithm is used. Throughout this paper, we analyze the tangent version of LSS, but the conclusions are also valid for the adjoint version since both algorithms give the same result up to round-off errors.

\section{How windowing mitigates the effect of time dilation}\label{proof}

The method introduced in this paper is similar to the original Least
Squares Shadowing method, except for two major differences.  The first
difference is the use of a smooth windowing function $\w$ satisfying
\[ \w(0)=\w(1)=0 \quad\mbox{and}\quad \int_0^1 \w(r) dr = 1; \]
which averages to 1 and tapers off to 0 at both ends of the
interval $[0,1]$.  The second difference is the lack of a time dilation
term in Equation (\ref{kkttan}) (Equation (\ref{kktadj}) for the adjoint version). Removal of
the time dilation term can simplify the implementation of a chaotic
sensitivity analysis capability for many solvers.\\

To understand the removal of the time dilation term, we should first understand
the original Least Squares Shadowing formulation \cite{explosion}, and why it
has a time dilation term.
The Least Squares Shadowing method approximates the derivative of a long
time average via an exchange of limit and derivative :
\begin{equation} \label{limitderiv}
    \frac{d\langle J\rangle}{ds}
    := \frac{d}{ds} \lim_{T\to\infty} \frac1T \int_0^T J(u(t;s),s) dt
     = \lim_{T\to\infty} \frac{d}{ds} \frac1T \int_0^T J(u(t;s),s) dt
\end{equation}
As a sufficient condition for this exchange, we require the total
derivative of $J$ with respect to $s$,
\[ \frac{d}{ds}J(u(t;s),s) = J_u \frac{\partial u}{\partial s} + J_s \]
to be uniformly continuous.
For this to hold, we ask for $J$ to have uniformly
continuous partial
derivatives, and $\frac{\partial u}{\partial s}$ to be uniformly
bounded and continuous. Among these conditions, it is most difficult to make
$\frac{\partial u}{\partial s}$ uniformly bounded, when $u$ is governed
by a family of chaotic dynamical systems.  $\frac{\partial u}{\partial
s}$ describes the difference between two solutions that satisfy
the governing equations with infinitesimally different parameters $s$ and $s+\delta s$.
It must satisfy the linearized governing equation
\begin{equation} \label{linearizednoeta}
\frac{\partial \vv}{\partial t} = f_u \vv + f_s .
\end{equation}
It is well known, as the
``butterfly effect'', that solutions to this linearized equation
can grow exponentially as $t$ increases. Under some assumptions on the dynamical system, the shadowing lemma (\cite{hyperbolic}, \cite{quasi}) ensures the existence of a non-exponentially growing solution. Nevertheless, this solution is potentially linearly growing and thus not uniformly bounded. For a uniformly bounded solution to exist, Equation (\ref{linearizednoeta}) has to be slightly modified to : 
\begin{equation} \label{linearizedeta}
\frac{\partial v}{\partial t} = f_u v + f_s + \eta f.
\end{equation}
for some uniformly bounded scalar function of time $\eta$, \cite{continuous}.
The solution to this equation, $v$, describes the difference between two
solutions that satisfy the governing equation with both infinitesimally
different parameter $s$ and infinitesimally different rate of time
marching.  Specifically, if $u(t)$ satisfies the equation
\[ \frac{du}{dt} = f(u,s), \]
then $u(t) + \epsilon v(t)$ satisfies the equation
\[ \frac{d(u + \epsilon v)}{(1 + \epsilon\eta) dt} = f(u+\epsilon v,s+\epsilon) \]
The relative difference between the rates of time marching is quantified by $\eta$.
The shadowing lemma ensures that there exists a pair of $v$ and $\eta$
that are both uniformly bounded in time and satisfy Equation
(\ref{linearizedeta}). The original Least Squares Shadowing method
approximates this uniformly bounded pair by minimizing the norm of
$v$ and $\eta$ under the constraint of Equation (\ref{linearizedeta}).
The desired derivative is then computed by modifying Equation
(\ref{limitderiv}) into :
\begin{equation} \label{modified}
    \frac{d\langle J\rangle}{ds}
= \langle J_u v\rangle + \langle J_s\rangle + \langle \eta J\rangle
- \langle\eta\rangle\langle J\rangle,
\end{equation}
to account for the time dilation effect of $\eta$.\\

The new method in this paper avoids time dilation by only computing an approximation of $v$ instead of the uniformly bounded pair $(v,\eta)$. This new quantity denoted $\vv_{\check{\tau}}$ is equal to :
\begin{equation} \label{vv0}
    \vv_{\check{\tau}}(t) = v(t) - f(t) \int_{\check{\tau}}^t \eta(r) dr
\end{equation}
for some ${\check{\tau}}$, where $f(t)$ denotes $f(u(t;s),s)=\frac{du}{dt}$. The subscript ${\check{\tau}}$ indicates that $\vv_{\check{\tau}}$ is "anchored" to $v$ at time ${\check{\tau}}$ since $\vv_{\check{\tau}}({\check{\tau}})=v({\check{\tau}})$. As we will see later on in section \ref{s:numericalresults} (figure \ref{lineargrowth}), the norm of $\vv_{\check{\tau}}$ has a V-shape and the "bottom of the V" is reached around $t={\check{\tau}}$ which explains the choice of the notation ${\check{\tau}}$. In general, $\vv_{\check{\tau}}$ has linear growth with respect to time due to the second term in equation (\ref{vv0}).  Time differentiating Equation (\ref{vv0}), using (\ref{linearizedeta}) and the fact that $\frac{df}{dt} = f_u \frac{du}{dt} = f_u f$, we recover :
\[ \frac{d \vv_{\check{\tau}}}{d t} = \frac{d v}{d t}
- \frac{d f}{d t} \int_{\check{\tau}}^t \eta dr
- \eta f = f_u \vv_{\check{\tau}} + f_s, \]
meaning that $\vv_{\check{\tau}}$ satisfies the simple no time-dilated Equation (\ref{linearizednoeta}).\\

Since this $\vv_{\check{\tau}}(t)$ defined by Equation (\ref{vv0}) is not uniformly bounded it
cannot be directly used to compute the desired derivative
by commuting the limit and derivative as in Equation (\ref{limitderiv}).
Instead, we introduce an approximation that involves a window
function $\w : [0,1]\to\R$ to mitigate the linear growth of $\vv_{\check{\tau}}$.  This approximation is 
\begin{equation}
    \frac{d\langle J\rangle}{ds}
    \approx \frac1T \int_0^T \w(\tfrac{t}{T}) (J_u
    \vv_{\check{\tau}} + J_s) dt
\end{equation}
for any ${\check{\tau}}\in[0,T]$.
The validity of this approximation is established through the following
theorem.

\begin{thm}
If the following are true:
\begin{enumerate}
        \item $\w$ is continuously differentiable,
        \item $\w(0) = \w(1) = 0$, and
        \item $\int_0^1 \w(r)dr = 1$,
        \item ${\check{\tau}}$ is a function of $T$ satisfying
            $0\le {\check{\tau}}\le T$ and
            $\displaystyle\lim_{T\to\infty}\frac{\check{\tau}}{T}$
            exists.
\end{enumerate}
then,
\begin{equation} \label{limitderivwin}
    \frac{d}{ds} \lim_{T\to\infty} \frac1T \int_0^T J(u(t;s),s) dt
    = \lim_{T\to\infty} \frac1T \int_0^T \w(\tfrac{t}{T}) (J_u
    \vv_{\check{\tau}} + J_s) dt\;.
\end{equation}
\end{thm}

This equality is nontrivial.  To prove that it is true,
we first define the following.
\begin{defi}
For any continuous function $\w:[0,1]\to \R$, the mean of the window is
\[ \overline{\w} := \int_0^1 \w(r)dr; \]
the infinitely-long windowed time average of a signal $x(t)$ is
\[ \langle x\rangle_{\w} := \lim_{T\to\infty} \frac1T \int_0^T
\w(\tfrac{t}{T}) x(t) dt, \]
\end{defi}

A special case of the window function is $\w\equiv 1$ called \textit{square} window. The mean of
this window is 1; the infinitely-long windowed time average of a signal
$x(t)$ is simply its ergodic average, which we already denoted as :
\[ \langle x\rangle := \lim_{T\to\infty} \frac1T \int_0^T x(t) dt. \]

\begin{lem} \label{lem:win}
If $ x(t) $ is bounded and $\langle x\rangle$ exists, then $\langle x\rangle_{\w}$ also exists for any
continuous $\w$, and the following equality is true:
\[ \langle x\rangle_{\w} = \overline\w\; \langle x\rangle. \]
\end{lem}
The proof is given in \ref{appendix1}.\\

Note that Lemma \ref{lem:win} does not apply to the windowed average on the right
hand side of Equation (\ref{limitderivwin}) in Theorem 1.  This is because
the lemma requires $x(t)$ to be independent of the
averaging length $T$.  But in Theorem 1,
${\check{\tau}}$ and thus $J_u \vv_{\check{\tau}} + J_s$ depends on $T$.
To apply Lemma \ref{lem:win}, we must first decompose $\vv_{\check{\tau}}$
according to Equation (\ref{vv0}).  Lemma \ref{lem:win} can then be
applied to all but one of the components:
\begin{equation}\label{decomposevv}\begin{split}
  & \lim_{T\to\infty} \frac1T \int_0^T \w(\tfrac{t}{T}) (J_u
    \vv_{\check{\tau}} + J_s) dt \\
= \;& \langle J_u v\rangle_{\w}
    + \langle J_s\rangle_{\w}
    - \lim_{T\to\infty} \frac1T\int_0^T \w(\tfrac{t}{T}) J_u f(t)
    ({\textstyle\int_{\check{\tau}}^t \eta(r) dr}) dt \\
= \;& \langle J_u v\rangle
    + \langle J_s\rangle
    - \lim_{T\to\infty} \frac1T\int_0^T \w(\tfrac{t}{T}) 
    \Big({\textstyle\int_{\check{\tau}}^t \eta(r) dr}\Big) \frac{dJ}{dt}dt
\end{split}\end{equation}
if $\overline{\w} = 1$. Here we used the fact that :
\[ \frac{dJ}{dt}
 = J_u \frac{du}{dt} + J_s \cancelto{0}{\frac{ds}{dt}}
 = J_u f \]
We then apply integration by parts to the remaining windowed average,
and use the assumption that $\w(0)=\w(1)=0$ in Theorem 1, to obtain :
\begin{equation} \label{ibp} \begin{split}
    - & \lim_{T\to\infty} \frac1T\int_0^T \w(\tfrac{t}{T})
    \Big({\textstyle\int_{\check{\tau}}^t \eta(r) dr}\Big) \frac{dJ}{dt}dt \\
=\; & \lim_{T\to\infty} \frac1T\int_0^T \w(\tfrac{t}{T}) \eta(t) J(t)\, dt
    + \frac1T\int_0^T \frac{d}{dt}\w(\tfrac{t}{T})
    \Big({\textstyle\int_{\check{\tau}}^t \eta(r) dr}\Big) J(t)\, dt \\
=\; & \langle \eta J\rangle
    + \lim_{T\to\infty} \frac1T \int_0^T \w'(\tfrac{t}{T})
    \tfrac1T \Big({\textstyle\int_{\check{\tau}}^t \eta(r) dr}\Big) J(t)\, dt
\end{split}\end{equation}
Here $\w'$ is the derivative of the window function $\w$.
By substituting Equation (\ref{ibp}) into Equation (\ref{decomposevv}),
then comparing with Equation (\ref{limitderiv}),
we see that we can prove Theorem 1 by proving the equality
\begin{equation}\label{final}
    \lim_{T\to\infty} \frac1T \int_0^T \w'(\tfrac{t}{T})\tfrac1T 
    \Big({\textstyle\int_{\check{\tau}}^t \eta(r) dr}\Big) J(t)\, dt
=- \langle\eta\rangle\langle J\rangle \end{equation}
We now establish this equality, thereby proving Theorem 1, through
two lemmas:

\begin{lem}\label{lemma3}
If $\eta$ is bounded and $\langle\eta\rangle$ exists, then
\[ \lim_{T\to\infty}\left(\sup_{{\check{\tau}},t\in[0,T]}
        \left(\tfrac{1}{T}\left(\textstyle\int_{\check{\tau}}^t \eta(r) dr\right) - 
 \langle\eta\rangle \frac{t - {\check{\tau}}}{T}\right)\right) = 0. \]
\end{lem}
The proof of this lemma is given in \ref{appendix2}. This lemma establishes the equality that
\[ \lim_{T\to\infty} \frac1T \int_0^T \w'(\tfrac{t}{T})
    \tfrac1T \Big({\textstyle\int_{\check{\tau}}^t \eta(r) dr}\Big) J(t)\, dt
 = \langle\eta\rangle
 \lim_{T\to\infty} \frac1T \int_0^T \w'(\tfrac{t}{T})
    \frac{t - {\check{\tau}}}{T} J(t)\, dt \]
The remaining task in proving Equation (\ref{final}) is achieved by the
following lemma:
\begin{lem}
    If $J$ is bounded, $\w\in C^1[0,1]$, $\w(0)=\w(1)=0$, and
    $\displaystyle\lim_{T\to\infty}\frac{\check{\tau}}{T}$ exists, then
\[ \lim_{T\to\infty} \frac1T \int_0^T \w'(\tfrac{t}{T})
    \frac{t - {\check{\tau}}}{T} J(t)\, dt 
=- \overline{\w}\langle J\rangle \]
\end{lem}
\begin{pf}
    Let $\check{\tau} = \displaystyle\lim_{T\to\infty}\frac{\check{\tau}}{T}$.
    Because both $J$ and $\w'$ are bounded ($\w$ is continuously
    differentiable in a closed interval),
\[ \lim_{T\to\infty} \frac1T \int_0^T \w'(\tfrac{t}{T})
    \frac{t - {\check{\tau}}}{T} J(t)\, dt 
   = \lim_{T\to\infty} \frac1T \int_0^T \w'(\tfrac{t}{T})
   \left(\tfrac{t}{T} - \check{\tau}\right) J(t)\, dt  \]
    Define $ \w_{\check{\tau}}(r) = \w'(r) (r - \check{\tau}), $
    then Lemma 2 can turn the equality above into,
\[ \lim_{T\to\infty} \frac1T \int_0^T \w'(\tfrac{t}{T})
    \frac{t - {\check{\tau}}}{T} J(t)\, dt 
= \langle J\rangle_{\w_{\check{\tau}}} = \overline{\w_{\check{\tau}}}\langle
J\rangle ,\]
in which
    \[ \overline{\w_{\check{\tau}}} := \int_0^1 \w'(r) (r - \check{\tau}) dr
        = \int_0^1 \frac{d\w}{dr} r\; dr
        - \cancelto{0}{\int_0^1 \frac{d\w}{dr} \check{\tau}\; dr}
        = -\int_0^1 w(r) dr = -\overline{\w}
    \]
    \qed
\end{pf}
In our case, not only $\displaystyle\lim_{T\to\infty}\frac{\check{\tau}}{T}$ exist but we even have $\displaystyle\lim_{T\to\infty}\frac{\check{\tau}}{T} \to \frac{1}{2}$ as shown in \ref{appendix3}. This result comes from the fact that the computed $\vv_{\check{\tau}}$ minimizes $\int_0^T |\vv_{\check{\tau}}(t)|^2dt$. \\

Lemma 3-4 combines to prove Equation (\ref{final}), which
combines with (\ref{decomposevv}) and (\ref{ibp}) to :
\begin{equation}\label{vv}
  \lim_{T\to\infty} \frac1T \int_0^T \w(\tfrac{t}{T}) (J_u
    \vv_{\check{\tau}} + J_s) dt
= \langle J_u v\rangle + \langle J_s\rangle + 
\langle \eta J\rangle - \langle\eta\rangle\langle J\rangle\;.
\end{equation}
This, together with Equation (\ref{modified}),
derived in previous literature (\cite{explosion}, \cite{continuous}), proves that the desired
derivative $\frac{d\langle J\rangle}{ds}$ can be computed via a windowed
average of $J_u \vv_{\check{\tau}} + J_s$, where $\vv$ is a solution to
the tangent equation without time dilation.

\section{Numerical results}
\label{s:numericalresults}
Now, we are going to apply simplified LSS to a test case dynamical system known as \textit{Lorenz 63}. It is a $3$-dimensional autonomous differential equation parameterized by $\sigma$, $\beta$, $\rho$ and the dynamics happen to be chaotic when the parameters belong to a certain range of values. The governing equations are the following :
 $$\left \{ \begin{array}{ll}
\frac{dx}{dt}=\sigma(y-x)\\
\frac{dy}{dt}=x(\rho-z)-y\\
\frac{dz}{dt}=xy-\beta z\\
\end{array}\right. $$

Edward Lorenz introduced them in 1963 to model the atmospheric convection. The quantity of interest in this example is $\langle z \rangle = \lim_{T \to \infty}\frac{1}{T}\int_{t=0}^{T}z(t) dt$,  the time average of the component $z$ and the algorithm is used to compute $\frac{d\langle z \rangle}{d\rho}$, its derivative  with respect to $\rho$. In this case, the design paramater $s$ is equal to $\rho$ and $J(u,\rho)=z$ where $u=(x,y,z)$. While a variety of values for $\rho$ will be tested, the other two parameters are set to $\sigma=10$ and $\beta=\frac{8}{3}$. It has been shown that $\frac{d\langle z \rangle}{d\rho}$ is approximately equal to $1$ for a wide range of values of $\rho$ (our numerical applications will stay in this range)\cite{butterfly}.\\

\begin{figure}[htb]
\begin{center}
\includegraphics[scale=0.6]{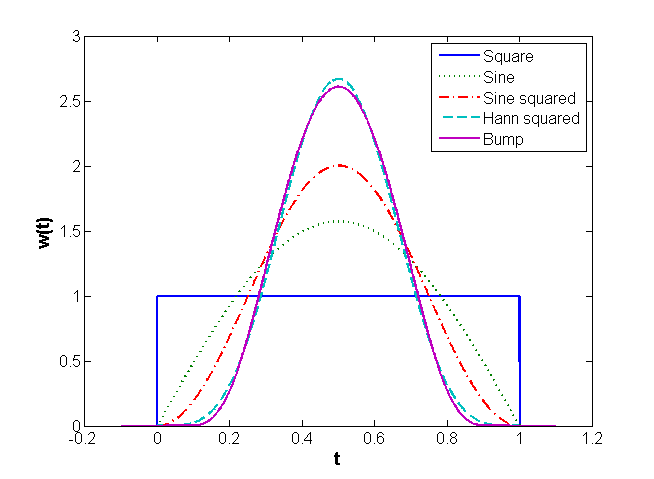}
\caption{Shape of the different windowing functions}\label{windowshape}
\end{center}
\end{figure}

The chosen windowing functions for this test case are (see figure \ref{windowshape}) :
\begin{itemize}
\item Square window : 
$$
\w(t)=1, \quad t\in [0,1]
$$
\item Sine window:
$$
\w(t)=\sin(\pi t), \quad t\in [0,1]
$$
\item Sine squared window :
$$
\w(t)=\sin^2(\pi t), \quad t\in [0,1]
$$
\item Sine to the power four also known as Hann-squared window :
$$
\w(t)=\sin^4(\pi t), \quad t\in [0,1]
$$
\item Bump window :
$$
\w(t)=\exp(\frac{-1}{t-t^2}), \quad t\in [0,1]
$$
defined up to normalization constants.
\end{itemize}
All windows respect the criteria of theorem 1 with the exception of the square window which is equivalent to not using any window. The major difference between the last four functions is their Taylor expansion in the neighborhood of $t=0$ and $t=1$. Actually, we can expect the areas around the extremities to be the most "sensible" ones, leading the the biggest error terms in the approximation of $\frac{d\langle z \rangle}{d\rho}$ since $\vv_{\check{\tau}}$ is maximal in these areas. Thus, choosing a windowing function that collapses quite fast on the extremities is expected to be the more accurate. The simple sine has non zero derivatives on the extremities, the squared one has a zero derivative on both sides, the Hann-squared function has zero coefficients until the 3rd order derivative while the bump function has zero coefficients for all derivative orders. A deep and rigorous analysis of this intuitive idea has been carried out by Krakos et al. for the periodic case (periodic dynamical system)\cite{periodic}.\\

For all simulations, the equations were discretized using an implicit second-order Crank-Nicolson scheme\footnote{It is worth noting that both integration scheme for $u(t)$ and discretization of equations (\ref{kkttan}) should have the same accuracy order for stability reasons.} solved up to machine precision and a uniform timestep of $\Delta t=0.02$. The burn-in period for computing $u(t)$ was set to $10$ time units.\\

\begin{figure}[htb]
\begin{center}
\includegraphics[scale=0.6]{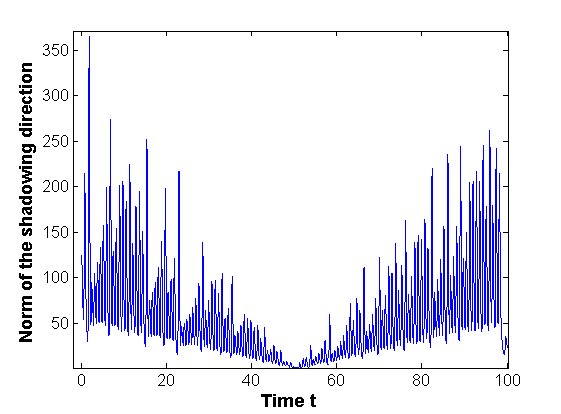}
\caption{Norm of $\vv_{\check{\tau}}$ with respect to $t$ for $\rho=28$}\label{lineargrowth}
\end{center}
\end{figure}

First, figure \ref{lineargrowth} shows the norm of the computed $\vv_{\check{\tau}}$ with respect to $t$ where $\rho=28$ and the time integration length is $T=100$. We can clearly notice the linear growth of the "envelope" of $\vv_{\check{\tau}}$ as $t$ increases or decreases as predicted in equation (\ref{vv0}) :
$$\vv_{\check{\tau}}(t)=v(t) - f(t) \int_{\check{\tau}}^t \eta(r) dr$$
Furthermore, $\vv_{\check{\tau}}$ is minimal around $t=\frac{T}{2}=50$ which is in agreement with the following result : 
$$\displaystyle\lim_{T\to\infty}\frac{\check{\tau}}{T} \to \frac{1}{2}$$

\begin{figure}[htb]
\begin{center}
\includegraphics[scale=0.65]{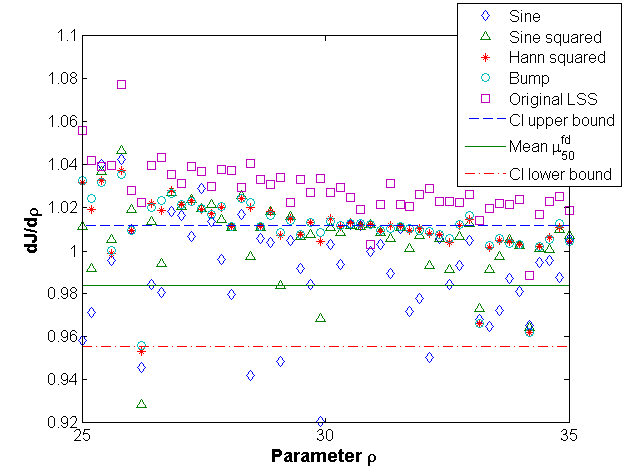}
\caption{$\frac{d\langle z \rangle}{d\rho}$ using different windowing functions and original LSS for integration time $T=50$}\label{rhospan}
\end{center}
\end{figure}

\begin{figure}
\centering
\includegraphics[width=10cm]{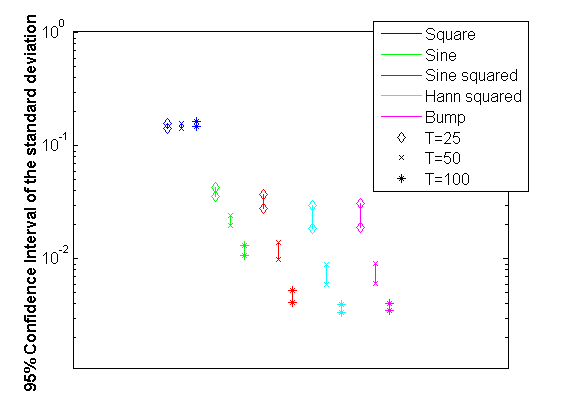}\\
\centering
\includegraphics[width=10.3cm]{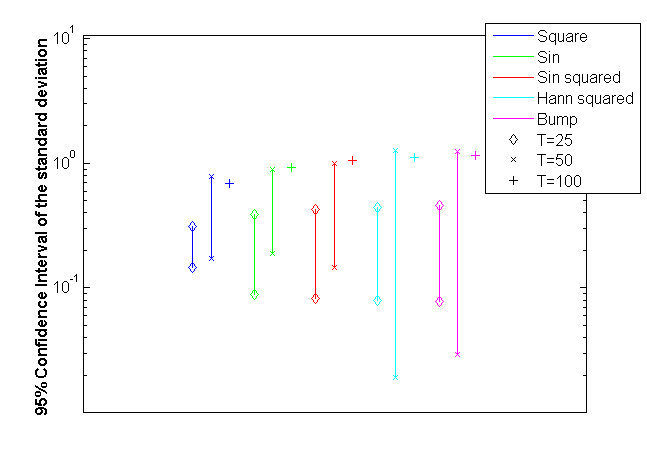}
\caption{Standard deviation confidence intervals for $\rho=28$ (upper graph) and $\rho=50$ (lower graph) using different windowing functions}\label{varplot}
\end{figure}

Then, we have computed $\frac{d\langle z \rangle}{d\rho}$ for $\rho \in [25,35]$ and $T=50$ using different windowing functions as well as the original Least Squared Shadowing algorithm (figure \ref{rhospan}). In this range of the parameter $\rho$, the dynamical system is known to be quasi-hyperbolic and the analysis we carried out in the previous section remains valid.\\
In order to have a rough estimation of $\frac{d\langle z \rangle}{d\rho}$ for this range of $\rho$, we computed $\langle z \rangle (\rho=25)$, $\langle z \rangle (\rho=35)$ for $50$ different initial conditions and constituted the finite difference approximation :
\begin{align}
\frac{d\langle z \rangle}{d\rho}\approx \frac{\langle z \rangle (\rho=35)-\langle z \rangle (\rho=25)}{35-25}
\end{align}  
Then, we represented in figure \ref{rhospan} the mean of these estimations as well as their $3$-sigma confidence interval:
\begin{align}
\textrm{CI}=\big[\mu_{50}^{\textrm{fd}} \pm 3\frac{\sigma_{50}^{\textrm{fd}}}{\sqrt{50}}\big]
\end{align}
where $\mu_{50}^{\textrm{fd}}$ is the average of the $50$ samples of $\frac{d\langle z \rangle}{d\rho}$ computed using finite difference and $\sigma_{50}^{\textrm{fd}}$ is their estimated standard deviation. After comparing these bounds with the results obtained with our algorithms, we notice that original LSS as well as the windowed algorithms give correct and very similar estimations of $\frac{d\langle z \rangle}{d\rho}$. Since the different versions of LSS have been tried on a single trajectory and for a finite integration time $T=50$, there is no reason for the computed derivatives to be unbiased. Furthermore, all windows do not seem to have the same performance: the results obtained with the Hann-squared and the bump window are smoother and more self-consistent than the ones coming from the sine squared and much more smoother than the ones given by the simple sine window. In order to point out this phenomenon, a complementary analysis has been done: for $\rho=28$, $\rho=50$ and integration time lengths of $T=25$, $T=50$ and $T=100$, $2000$ simulations (each one with a different initial condition) were run to compute a $95\%$ confidence interval of the standard deviation $\sigma$ of $\frac{d\langle z\rangle}{d\rho}$ for each one of the windows (see figure \ref{varplot}). The $95\%$ confidence intervals were approximated as follows :
\begin{align}
\textrm{CI}=\Bigg[\Bigg(\sigma^2_{2000} \pm 1.96\sqrt{\frac{\sigma_{2000}^4}{2000}\big((\kappa_{2000}-1)+\frac{2}{2000-1}\big)}\Bigg)^{\frac{1}{2}}\Bigg]
\end{align}
where $\sigma_{2000}$ and $\kappa_{2000}$ are, respectively, the estimates of the 2000 samples standard deviation and kurtosis.\\
The results confirm our previous remark: for a fixed $T$ and $\rho=28$, the Hann-squared and bump windows have the lowest standard deviations meaning they are more robust than the other windows. Among the valid windowing functions, the simple sine gives the worst results. Then, for a fixed window, we notice that increasing the time integration length decreases the standard deviation which is an intuitive result since when $T$ gets bigger the influence of the initial condition on the dynamical system fades out (due to ergodicity). Finally, for $\rho=50$, all windows give bad results (standard deviation of order 1) which can be explained by the fact that the dynamical system is no longer quasi hyperbolic for this value of $\rho$. The theory we have developed doesn't hold anymore and there is no reason for simplified LSS to converge to the true value of $\frac{d\langle z \rangle}{d\rho}$ (for $T=100$, the confidence intervals were so big and uninformative that we only represented the standard deviation estimator).

\section{Conclusion}
Simplified LSS is an improvement of original LSS in the sense that the time dilation factors have been removed from the formulation of the problem leading to a much simpler system of equations to solve. When running the original version of LSS, we are lead to solve a slightly more complex set of KKT equations than (\ref{kkttan}) which can be found in \cite{lsspde}. This set of equations is solved by inverting the Schur complement of the KKT matrix which contains the concatenation of smaller matrices of the form $f_{u_i} \times f_{u_i}^T + f(u_i)\times f(u_i)^T$ where $\{u_i\}$ is a finite discretization of a trajectory. The first term is a multiplication of two $n\times n$ matrices where $n$ is the phase space dimension. For PDEs, $n$ can be extremely large. However, since the finite stencil used to discretize the PDE is local, $f_{u_i}$ is very sparse, so is $f_{u_i} \times f_{u_i}^T$. The second term, which is totally due to the time dilation terms, is a $n\times 1$ by $1\times n$ multiplication which leads to a rank one but fully dense $n\times n$ matrix. Consequently, $f_{u_i} \times f_{u_i}^T + f(u_i)\times f(u_i)^T$ is dense as well and the Schur complement becomes hard, if not impossible to inverse for high-dimensional systems.\\
With windowing, the time dilation factors are no longer needed and the Schur complement is now only formed of $f_{u_i} \times f_{u_i}^T$ submatrices. We can consequently use sparse solvers and preconditioners such as incomplete LU factorization to solve the KKT system. Simplified LSS can be efficiently applied to very large, industrial-scale problems which were very hard if not impossible to tackle with original LSS \cite{fluidexample}. Additionally, simplified LSS is easier to implement.\\

The simplification has been made possible thanks to the introduction of a windowing function that mitigates the linear growth of the shadowing direction $\vv_{\check{\tau}}$. This function is required to be continuously differentiable, null at its extremities and with integral equal to $1$. Once picked in this class of functions, the major characteristic of the window is the way it decays at its extremities: a rapidly decaying function such as the bump window or the Hann-squared one will be less sensitive to transient behavior and give more robust results. Since all windowing functions lead to the same computational cost, our study strongly recommends the choice of a fastly decaying window. In the same way, increasing the time integration length decreases the influence of the initial condition making the algorithm more robust and accurate.   

\section*{Acknowledgements}
We acknowledge AFOSR Awards 14RT0138 under Dr. Fariba Fahroo and Dr. Jeanluc Cambrier, F11B-T06-0007 under Dr. Fariba Fahroo, and NASA Award NNH11ZEA001N under Dr. Harold Atkins.

\appendix
\section{Derivation of the adjoint version of \textit{simplified LSS}}\label{adjointderivation}
The adjoint version is directly obtained from the equations of the tangent version. We should have in mind that the objective function we want to compute is :
\begin{align}\label{objfun}
\frac{d\J}{ds}=\frac1T\int_0^T
            \w(\tfrac{t}{T}) \left(J_u \vv + J_s \right) dt
\end{align}
In order to obtain the adjoint formulation, we need to replace the first term $\w(\tfrac{t}{T}) J_u \vv$ which is based on the solution $\vv$ to the forward problem by a new term based on a backward solution $\hat{\vv}$. 
First, we reformulate the \textit{forward} set of equations (\ref{kkttan}):
\begin{equation} \label{kkttan2} \begin{cases}
            \frac{d\vv}{dt} - f_u \vv - f_s = 0 & \\
            \frac{d\ww}{dt} +f_u^T \ww - \vv = 0 &  \\
            \ww(0) = \ww(T) = 0             & 
\end{cases}\end{equation}
Then, we introduce the dual variables $\hat{\vv}$ and $\hat{\ww}$ for the first two equations and obtain the scalar equation :
\begin{align}
\int_0^T \hat{\vv}^T(\frac{d\vv}{dt} - f_u \vv - f_s) + \hat{\ww}^T(\frac{d\ww}{dt} +f_u^T \ww - \vv)dt = 0
\end{align}
where the superscript $T$ stands for the transpose operation. After integrating by parts and thanks to the third equation of (\ref{kkttan2}) :
\begin{align}
0&=\hat{\vv}^T\vv(T)-\hat{\vv}^T\vv(0)+\int_0^T \hat{\vv}^T(-f_u \vv - f_s) - \frac{d\hat{\vv}}{dt}^T\vv + \hat{\ww}^T(f_u^T \ww - \vv) -\frac{d\hat{\ww}}{dt}^T\ww dt\\
&=\vv^T\hat{\vv}(T)-\vv^T\hat{\vv}(0)+\int_0^T (-f_u^T\hat{\vv}-\frac{d\hat{\vv}}{dt}-\hat{\ww})^T\vv + (f_u \hat{\ww}-\frac{d\hat{\ww}}{dt})^T\ww - f_s^T \hat{\vv} dt \label{equaadj} 
\end{align}
Based on (\ref{objfun}), we impose the \textit{adjoint} (or \textit{backward}) set of equations :
\begin{equation} \begin{cases} \label{setequatadj}
            -f_u^T\hat{\vv}-\frac{d\hat{\vv}}{dt}-\hat{\ww} = \w(\frac{t}{T}) J_u^T & \\
            f_u \hat{\ww}-\frac{d\hat{\ww}}{dt} = 0 &  \\
            \hat{\vv}(0) = \hat{\vv}(T) = 0             & 
\end{cases}\end{equation}
After replacing the equalities of (\ref{setequatadj}) into (\ref{equaadj}), we get:
\begin{align}\label{adjobj2}
\frac{1}{T}\int_0^T (-f_u^T\hat{\vv}-\frac{d\hat{\vv}}{dt}-\hat{\ww})^T\vv dt=\frac{1}{T}\int_0^T \w(\tfrac{t}{T}) J_u\vv dt=\frac{1}{T}\int_0^T f_s^T \hat{\vv} dt
\end{align}
Combining (\ref{objfun}) and (\ref{adjobj2}), the new objetive function becomes:
\begin{align}
\frac{d\J}{ds}=\frac{1}{T}\int_0^T \left(f_s^T \hat{\vv}+ J_s \right) dt
\end{align}

\section{Lemma 2}\label{appendix1}
Given any $ \w \in C[0,1], \w\geq 0$, and bounded function $ x(t) $ such that 
$ \langle x\rangle = \lim_{T\to\infty}\frac1T \int_0^T  x(t) dt $ exists, there is the relation:

\begin{equation} \label{eqni}	
    \lim_{T\to\infty}\frac1T \int_0^T \w(\tfrac{t}T) x(t) dt
  = \overline{\w}\cdot \langle x\rangle
\end{equation}
\begin{pf}
	Given any $\epsilon>0$,  let $\epsilon_1 = \epsilon/2B_x$, here $B_x$ is the bound for $x(t)$.	
	$ \w $ is continuous  on a compact set, hence its range is also bounded, say, by $B_w$. 
	Moreover, $ \w $ is uniformly continuous, hence $\exists K$, 
	s.t. $\forall \tau_1,\tau_2 \in [0,1], \left| \tau_1 - \tau_2 \right| \leq 1/K$, 
	we have  $\left|  \w(\tau_1) - \w(\tau_2) \right| \leq \epsilon_1 $.
	Construct the partition $P=\{0,\frac{1}{K},\frac{2}{K}, ..., 1\}$. 
	Now $\forall i \in\{1,2,...,K\}$,
	let $ \w_i = K \cdot \int_{(i-1)/K}^{i/K} w(\tau)d\tau $, 
	then $ \overline{\w}=  \frac1K \sum_{i=1}^{K}  \w_i $.
	Let	$\xi_i(\tau) =  \w(\tau) - \w_i $, 
	$\tau \in I_i = \left[(i-1)/K, i/K \right]$. 
	
	$\forall i, \w $ is continuous on the segment $I_i$,
	which is a compact set,
	hence $ \w $ achieves its maximum and minimum value on this segment.
	Denote its maximum and minimum values on this segment by $ \w_{i,max} $ and $ \w_{i,min} $ respectively. 
	By selection of K, $ \w_{i,max} - \w_{i,min} \leq \epsilon_1$.
	By definition of $ \w_i $, $ \w_{i,min} \leq \w_i \leq \w_{i,max}$. 
	Hence we have $ \xi_i(\tau) \leq \epsilon_1, \tau \in I_i$.

	With above $ K $, let $\epsilon_2 = \epsilon / 2KB_w $. 
	Since $ \lim_{T\to\infty}\frac1T \int_0^T  x(t) dt $ exists, 
	$ \exists T_0 $, 
	such that  $\forall M \geq T_0 $ , $\langle x\rangle - \frac1M \int_0^M  x(t) dt \leq \epsilon_2$. Hence we have:
	\[ \left|  \langle x\rangle - \frac{1}{iM} \int_{0}^{iM}  x(t) dt \right|\leq \epsilon_2, \quad i = 1,2,3,...  \]
	the i-th inequality is equivalent to:
	\[ 
	-i \epsilon_2 \leq  i \langle x\rangle - \frac{1}{M} \int_{0}^{M}  x(t) dt...- \frac{1}{M} \int_{(i-1)M}^{iM}  x(t) dt \leq i \epsilon_2, \quad i = 1,2,3,...
	\]
	the (i+1)-th inequality is equivalent to:
	\[ 
	-(i+1) \epsilon_2 \leq  -(i+1) \langle x\rangle + \frac{1}{M} \int_{0}^{M}  x(t) dt...+ \frac{1}{M} \int_{iM}^{(i+1)M}  x(t) dt \leq (i+1) \epsilon_2, \quad i = 0,1,2,...
	\]
	add above two inequalities together,
	\[ 
	- (2i+1)\epsilon_2 \leq  \langle x\rangle - \frac{1}{M} \int_{iM}^{(i+1)M}  x(t) dt \leq  (2i+1)\epsilon_2, \quad i = 0,1,2,...
	\]
	
	Now $\forall T \geq KT_0$, let $M = T/K$, then $ M\geq T_0$, 
	and the difference between the two sides in \ref{eqni}, when $ T $ is finite, could be represented by:
	\begin{align*}
	\overline{\w}\cdot \langle x\rangle	&- \frac1T \int_0^T \w(\tfrac{t}T) x(t) dt \\
	&=  \frac1K \sum_{i=1}^{K}  \left[ 
	\w_i \cdot \langle x\rangle
	- \frac1{M}  \int_{(i-1)M}^{iM} (\w_i+\xi_i(\tfrac{t}T))  x(t) dt  
	\right]  \\
	&= \frac1{K} \sum_{i=1}^{K}    \left[ 
	- \frac 1M \int_{(i-1)M}^{iM} \xi_i \left(\tfrac{t}{T}\right) x(t) dt  
	+ \w_i \cdot \left( \langle x\rangle - \frac1{M} \int_{(i-1)M}^{iM}  x(t) dt  \right)
	\right] 
	\end{align*}
	the first term is confined by:
	\[
	\left| \frac 1M \int_{(i-1)M}^{iM} \xi_i (\tfrac{t}{T}) x(t) dt  
	\right| 
	\leq
	\epsilon_1 \cdot \frac 1M\int_{(i-1)M}^{iM} \left| x(t) \right| dt  	 	
	\]
	the second term is confined by:
	\[
	\left| \w_i \cdot \left( \langle x\rangle - \frac1{M} \int_{(i-1)M}^{iM}  x(t) dt  \right) \right|
	\leq
	\w_i   \cdot (2i-1)\epsilon_2
	\]
	As a result,
	\begin{align*}
	\biggl| \overline{\w}\cdot \langle x\rangle	 &- \frac1T \int_0^T \w(\tfrac{t}T) x(t) dt \biggr| \\
	& \leq \frac1{K} \sum_{i=1}^{K}  \left[ 
	\epsilon_1 \cdot \frac 1M\int_{(i-1)M}^{iM} | x(t) | dt
	+ 
	\w_i   \cdot (2i-1) \epsilon_2 
	\right] \\
	& \leq \epsilon_1 \cdot B_x 
	+
	\epsilon_2 \cdot  B_w K 
	\leq \epsilon
	\end{align*} 
	
	Now we find $ T^* = KT_0$, s.t. $\forall T\geq T^*,$ $\left|  \frac1T \int_0^T \w(\tfrac{t}T) x(t) dt -  \overline{\w}\cdot \langle x\rangle\right|  \leq \epsilon$.
	\qed
\end{pf}

\section{Lemma \ref{lemma3}}\label{appendix2}
If $\eta$ is bounded and $\langle\eta\rangle$ exists, then
\[ \lim_{T\to\infty}\left(\sup_{{\check{\tau}},t\in[0,T]}
        \left(\tfrac{1}{T}\left(\textstyle\int_{\check{\tau}}^t \eta(r) dr\right) - 
 \langle\eta\rangle \frac{t - {\check{\tau}}}{T}\right)\right) = 0. \]
\begin{pf}
First, we write:
\begin{align}
\tfrac{1}{T}\left(\textstyle\int_{\check{\tau}}^t \eta(r) dr\right)-\langle\eta\rangle \frac{t-{\check{\tau}}}{T}=\tfrac{t-{\check{\tau}}}{T}\bigg(\frac{1}{t-{\check{\tau}}}\left(\textstyle\int_{\check{\tau}}^t \eta(r) dr\right)-\langle\eta\rangle\bigg)
\end{align}
Based on the two assumptions  $t,{\check{\tau}}\in[0,T]$ and $\eta\leq \|\eta\|^\infty<\infty$, we obtain the following bounds on the two terms of the product we have just derived:

\begin{align}
\Big| \frac{t-{\check{\tau}}}{T} \Big|\leq 1 \qquad \textrm{as well as} \qquad \Big|\frac{1}{t-{\check{\tau}}}\left(\textstyle\int_{\check{\tau}}^t \eta(r) dr\right)-\langle\eta\rangle\Big| \leq 2\| \eta \|^\infty
\end{align}
Let us define $c=\max\big(1,2\|\eta\|^\infty\big)$. Based on the ergodicity of $\eta$, for any $\epsilon>0$ there is a constant $M>0$ such that for all $t_0$ and $T\geq M$ :
\begin{align}\label{uniformergod}
\Big| \frac{1}{T}\int_{r=t_0}^{t_0+T} \eta(r) dr  - \langle \eta \rangle \Big |< \epsilon
\end{align}

Consequently, for all $\epsilon>0$ and for all $T>\frac{cM}{\epsilon}$, we have:
\begin{itemize}
\item if $|t-{\check{\tau}}|<M$, then $\Big|\frac{t-{\check{\tau}}}{T}\Big|\leq \frac{\epsilon}{c}$ which means that :
\begin{align}
\Bigg|\tfrac{t-{\check{\tau}}}{T}\bigg(\frac{1}{t-{\check{\tau}}}\left(\textstyle\int_{\check{\tau}}^t \eta(r) dr\right)-\langle\eta\rangle\bigg)\Bigg| \leq \frac{\epsilon}{c}\Bigg|\frac{1}{t-{\check{\tau}}}\left(\textstyle\int_{\check{\tau}}^t \eta(r) dr\right)-\langle\eta\rangle\Bigg| \leq \epsilon
\end{align}
 \item if $|t-{\check{\tau}}|\geq M$, thanks to relation (\ref{uniformergod}):
 \begin{align}
 \Big|\frac{1}{t-{\check{\tau}}}\left(\textstyle\int_{\check{\tau}}^t \eta(r)    dr\right)-\langle\eta\rangle\Big|\leq \epsilon
 \end{align}
 which again implies :
\begin{align}
\Big|\tfrac{t-{\check{\tau}}}{T}\bigg(\frac{1}{t-{\check{\tau}}}\left(\textstyle\int_{\check{\tau}}^t \eta(r) dr\right)-\langle\eta\rangle\bigg)\Big|\leq \epsilon
\end{align}
\end{itemize}
that concludes the proof.

\end{pf}

\section{Proof of convergence of $\lim_{T\to\infty}\frac{\check{\tau}}{T}$}\label{appendix3}
When minimizing the cost function $\int_0^T\|\check{v}_{\check{\tau}}(t)\|^2dt$, the choice of ${\check{\tau}}$ is arbitrary, thus :
\begin{align}
\frac{d}{d{\check{\tau}}}\Big(\int_0^T\|\check{v}_{\check{\tau}}(t)\|^2dt\Big)=0
\end{align}
at the global minimum. We have :
\begin{align}
\int_0^T  \| \check{v}_{\check{\tau}}(t)\|^2 dt = \int_0^T \Big(\|v(t)\|^2+ \| f(t)\|^2\big(\int_{\check{\tau}}^t \eta(r)dr\big)^2-2v(t)^T f(t)\int_{\check{\tau}}^t \eta(r)dr \Big)dt
\end{align}
Thus :
\begin{align}
\frac{d}{d{\check{\tau}}}\Big(\int_0^T\|\check{v}_{\check{\tau}}(t)\|^2dt\Big)&=-2\int_0^T\|f(t)\|^2\eta({\check{\tau}})\big(\int_{\check{\tau}}^t \eta(r)dr\big)dt +2\int_0^Tv(t)^Tf(t)\eta({t})dt\\
&=2\eta({\check{\tau}})\Big(\int_0^Tv(t)^Tf(t)dt-\int_0^T\|f(t)\|^2\big(\int_{\check{\tau}}^t \eta(r)dr\big)dt\Big)
\end{align}
It can be shown that the ${\check{\tau}}$ such that $\eta({\check{\tau}})=0$ if any, don't correspond to a global minimum. We consequently should have :
\begin{align}
\label{equality}
\int_0^Tv(t)^Tf(t)dt-\int_0^T\|f(t)\|^2\big(\int_{\check{\tau}}^t \eta(r)dr\big)dt=0
\end{align}
On one side :
\begin{align}
\Big|\int_0^Tv(t)^Tf(t)dt\Big|&\leq \int_0^T \|v\|^\infty \|f\|^\infty dt\\
 &\leq T \|v\|^\infty \|f\|^\infty
\end{align}
On the other, $\int_{\check{\tau}}^t \eta(r)dr \sim (t-{\check{\tau}})\langle \eta \rangle$ which implies :
\begin{align}
\Big| \int_0^T\|f(t)\|^2\big(\int_{\check{\tau}}^t \eta(r)dr\big)dt \Big|\sim \Big|\int_0^T\|f(t)\|^2(t-{\check{\tau}})\langle \eta \rangle dt\Big|
\end{align}
after integrating by parts :
\begin{align}
 \Big|\int_0^T\|f(t)\|^2(t-{\check{\tau}})\langle \eta \rangle dt\Big| &= |\langle \eta \rangle |\Big| (T-{\check{\tau}})\int_0^T\|f(s)\|^2ds - \int_0^T\int_0^t\|f(s)\|^2dsdt\Big|\\
&\sim  |\langle \eta \rangle | \Big| T(T-{\check{\tau}})\langle \|f\|^2\rangle - \frac{T^2}{2}\langle \|f\|^2\rangle\Big|\\
&\sim \Big| \langle \eta \rangle \langle \|f\|^2\rangle T^2(\frac{1}{2}-\frac{\check{\tau}}{T})\Big|
\end{align}
which means that  :
\begin{align}
\Big| \int_0^T\|f(t)\|^2\big(\int_{\check{\tau}}^t \eta(r)dr\big)dt \Big| \sim | \langle \eta \rangle \langle \|f\|^2\rangle | \times T^2|(\frac{1}{2}-\frac{\check{\tau}}{T})|
\end{align}
This relation shows that if $\frac{\check{\tau}}{T}\not\to\frac{1}{2}$ as $T\to \infty$, the second term grows as $T^2$ while the first term grows, at most, as $T$. This contradicts equation  (\ref{equality}) and concludes the proof.

\end{document}